\documentclass[11pt]{article}

\usepackage{amssymb}
\usepackage{graphicx}
\usepackage[mathscr]{eucal}
\usepackage[all]{xy}
\usepackage{amsmath}
\usepackage{pgfplots}
\usepackage{tikz}
\usetikzlibrary{calc}

\begin{document}

\centerline{ \bf \large ON THE COMPATIBILITY OF B\'{E}ZOUT COEFFICIENTS}
\centerline{\bf \large BETWEEN PYTHAGOREAN PAIRS}
\centerline{\bf \large UNDER UNIMODULAR TRANSFORMATIONS}

\vskip 0.2in
\centerline{Cherng-tiao Perng and Maila Brucal-Hallare}
\vskip 0.2in
\noindent {\small \textsc{Abstract.} In a recent preprint, Gullerud and Walker [2] proved a theorem and made a conjecture about the correctness of efficiently generating B\'{e}zout trees for Pythagorean pairs. In this note, we give a simple proof of their theorem, confirm that their conjecture is true, and furthermore we give a generalization.}\\

\noindent {\bf 1. Introduction} \\

\noindent The integers triple $(x,y,z)$ is called a Pythagorean triple if $x^2+y^2=z^2$. It is called primitive if they are relatively prime. It is well known that all positive primitive Pythagorean triples $(x,y,z)$ with $y$ even can be written as $$x=m^2-n^2,y=2mn,z=m^2+n^2,$$ for some relative prime integers $m$ and $n$ such that $m>n>0$ [4]. Following the authors of [2], we call such $(m,n)$ a \emph{Pythagorean pair}. Given $(m,n)$, it is clear that $(n,m)$ and $(m,-n)$ also generate Pythagorean triples; such pairs are called \emph{associated pairs} of $(m,n).$ Note that if $(m,n)$ is a Pythagorean pair, then $f(m,n):=(2m+n,m)$ (where $f$ is defined on ${\mathbb Z}\times{\mathbb Z}$) is another Pythagorean pair. Similarly, $f(n,m)=(2n+m,n)$ and $f(m,-n)=(2m-n,m)$ are also Pythagorean pairs. Define now a trinary tree generated by $(m,n)$ as follows:

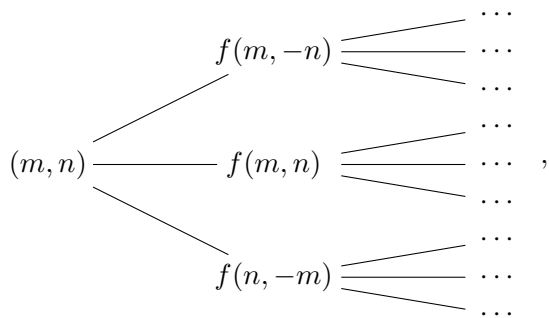
\def\a{2}
\def\b{0.5}
\def\c{0.65}
\def\d{0.7}
\def\e{3}

\begin{center}
\begin{tikzpicture}

\draw (0,0) node {$(m,n)$};

\draw (\e,3*\b) node {$f(m,-n)$};
\draw (\e,0) node {$f(m,n)$};
\draw (\e,-3*\b) node {$f(n,-m)$};

\draw (2*\e,4*\b) node {$\cdots$};
\draw (2*\e,3*\b) node {$\cdots$};
\draw (2*\e,2*\b) node {$\cdots$};
\draw (2*\e,\b) node {$\cdots$};
\draw (2*\e,0) node {$\cdots$}; \draw (2.2*\e,0) node {$,$};
\draw (2*\e,-\b) node {$\cdots$};
\draw (2*\e,-2*\b) node {$\cdots$};
\draw (2*\e,-3*\b) node {$\cdots$};
\draw (2*\e,-4*\b) node {$\cdots$};

\draw (0+0.2*\e,0+0.2*3*\b) -- (\e-0.2*\e,3*\b-0.2*3*\b);
\draw (0+0.2*\e,0) -- (\e-0.25*\e,0);
\draw (0+0.2*\e,0-0.2*3*\b) -- (\e-0.2*\e,-3*\b+0.2*3*\b);

\draw (\e+0.3*\e,3*\b+0.3*\b) -- (2*\e-0.15*\e,4*\b-0.15*\b);
\draw (\e+0.3*\e,3*\b) -- (2*\e-0.15*\e,3*\b);
\draw (\e+0.3*\e,3*\b-0.3*\b) -- (2*\e-0.15*\e,2*\b+0.15*\b);

\draw (\e+0.3*\e,0+0.3*\b) -- (2*\e-0.15*\e,\b-0.15*\b);
\draw (\e+0.3*\e,0) -- (2*\e-0.15*\e,0);
\draw (\e+0.3*\e,0-0.3*\b) -- (2*\e-0.15*\e,-\b+0.15*\b);

\draw (\e+0.3*\e,-3*\b+0.3*\b) -- (2*\e-0.15*\e,-2*\b-0.15*\b);
\draw (\e+0.3*\e,-3*\b) -- (2*\e-0.15*\e,-3*\b);
\draw (\e+0.3*\e,-3*\b-0.3*\b) -- (2*\e-0.15*\e,-4*\b+0.15*\b);

\end{tikzpicture}
\end{center}

\noindent where recursively, each node on a given level produces three nodes on a next level by applying the functions $f_1(m,n)=f(m,-n),f_2(m,n)=f(m,n)$ and $f_3(m,n)=f(n,m)$ and so on. Randall and Saunders [5] proved that the trinary tree produced from $(3,1)$ contains all pairs of relatively prime odd integers. Similarly, the trinary tree produced from $(2,1)$ contains all pairs of relatively prime integers of opposite parity. Thus these together generate all relatively prime Pythagorean pairs $(m,n)$ with $m>n>0$.\\

\noindent We call $(r,s)$ the B\'{e}zout coefficients associated with $(m,n)$ if $(r,s)$ is obtained from the standard division algorithm so that $rm+sn=\gcd(m,n)$. For comparison, for an input of $(m,n)$ in the Matlab \verb"gcd" function $$[G,U,V]=\gcd(m,n),$$ the output will be $G=\gcd(m,n)$ in the usual notation, and $U=r,V=s$ are the B\'{e}zout coefficients. In an attempt to efficiently generate the B\'{e}zout coefficients for Pythagorean pairs, Guillerud and Walker introduced the notion of \textbf{B\'{e}zout tree of $(m,n)$ generated by $(u,v)$}, which is defined by $$g(u,v)=(v,u-2v),$$
                                          $$g(v,u)=(u,v-2u)~{\rm and~}$$
                                          $$g(u,-v)=(-v,u+2v),$$
and the tree is arranged in the analogous format as in the tree starting with $(m,n).$ Guillerud and Walker proved the following result, for which we offer a simple argument.\\

\noindent {\bf Theorem 1.1.} (cf. Theorem 1.2 of [2]) Let $(m,n)$ be a Pythagorean pair with $m>n$ with associated pairs $(n,m)$ and $(m,-n)$. Let $f$ and $g$ be as defined above, and let $mu+nv=1$ for some $u,v\in {\mathbb Z}$. Then $g(u,v),g(v,u)$ and $g(u,-v)$ respectively yield the necessary coefficients $u',v'$ such that $$(2m+n)u'+mv'=1,$$ $$(2n+m)u'+nv'=1,~{\rm and~}$$ $$(2m-n)u'+mv'=1$$ respectively.\\

\noindent {\it Proof.} In terms of matrices, we have $$f(m,n)=\left[\begin{array}{c}2m+n\\ m\end{array}\right]=\left[\begin{array}{cc}2&1\\ 1&0\end{array}\right]\left[\begin{array}{c}m\\ n\end{array}\right].$$ If we let $A=\left[\begin{array}{cc}2&1\\ 1&0\end{array}\right]$, then it is clear that $g(u,v)$ is given by $$g(u,v)=\left[\begin{array}{c}v\\ u-2v\end{array}\right]=\left[\begin{array}{cc}0&1\\ 1&-2\end{array}\right]\left[\begin{array}{c}u\\ v\end{array}\right]=(A^{-1})^T\left[\begin{array}{c}u\\ v\end{array}\right].$$ It follows that $\left[\begin{array}{c}u'\\ v'\end{array}\right]:=g(u,v)$ satisfies $$(2m+n)u'+mv'=[u',v']\left[\begin{array}{c}2m+n\\ m\end{array}\right]=g(u,v)^Tf(m,n)$$
$$=[u,v]A^{-1}A\left[\begin{array}{c}m\\ n\end{array}\right]=um+vn=1,$$ as required. The other two cases are handled in exactly the same way. $\Box$\\

\noindent Before stating the conjecture (and we call it Theorem 1.4 now), let's look at the following example (Example 1.3 of [2]), where on the left it is the trinary tree generated by $(3,1)$ up to a depth of $2$, and on the right, it is the B\'{e}zout tree of $(3,1)$ generated by $(0,1)$ up to the same depth. Note that the defining rule for the second tree is analogous to the first: one proceeds from one node at a given level to three nodes at the next level by applying the functions $g_1(u,v)=g(u,-v),g_2(u,v)=g(u,v)$ and $g_3(u,v)=g(v,u).$\\

\noindent {\bf Example 1.2.}\\

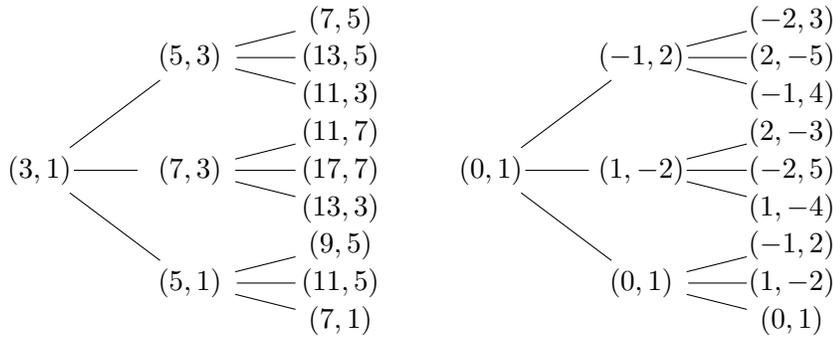
\noindent \begin{tikzpicture}

  \draw (0,0) node {$(3,1)$};
  \draw (0+1*\a,0) node {$(7,3)$};
  \draw (0+\a,3*\b) node {$(5,3)$};
  \draw (0+\a,-3*\b) node {$(5,1)$};
  \draw (0+2*\a,\b) node {$(11,7)$};
  \draw (0+2*\a,0) node {$(17,7)$};
  \draw (0+2*\a,-\b) node {$(13,3)$};
  \draw (0+2*\a,2*\b) node {$(11,3)$};
  \draw (0+2*\a,3*\b) node {$(13,5)$};
  \draw (0+2*\a,4*\b) node {$(7,5)$};
  \draw (0+2*\a,-2*\b) node {$(9,5)$};
  \draw (0+2*\a,-3*\b) node {$(11,5)$};
  \draw (0+2*\a,-4*\b) node {$(7,1)$};

  \draw (0+0.7*\c,0) -- (0+\a-\d,0);
  \draw (0+0.2*\a,0+0.2*3*\b) -- (0+0.8*\a,0+0.8*3*\b);
  \draw (0+0.2*\a,0-0.2*3*\b) -- (0+0.8*\a,0-0.8*3*\b);

  \draw (\a+0.3*\a,3*\b+0.3*\b) -- (2*\a-0.3*\a,4*\b-0.3*\b);
  \draw (\a+0.6\a,3*\b) -- (2*\a-0.3*\a,3*\b);
  \draw (\a+0.3*\a,3*\b-0.3*\b) -- (2*\a-0.3*\a,2*\b+0.3*\b);

  \draw (\a+0.3*\a,-3*\b-0.3*\b) -- (2*\a-0.3*\a,-4*\b+0.3*\b);
  \draw (\a+0.6\a,-3*\b) -- (2*\a-0.3*\a,-3*\b);
  \draw (\a+0.3*\a,-3*\b+0.3*\b) -- (2*\a-0.3*\a,-2*\b-0.3*\b);

  \draw (\a+0.3*\a,0+0.3*\b) -- (2*\a-0.3*\a,\b-0.3*\b);
  \draw (\a+0.3*\a,0) -- (2*\a-0.3*\a,0);
  \draw (\a+0.3*\a,0-0.3*\b) -- (2*\a-0.3*\a,-\b+0.3*\b);

  \draw (0+3*\a,0) node {$(0,1)$};
  \draw (0+1*\a+3*\a,0) node {$(1,-2)$};
  \draw (0+\a+3*\a,3*\b) node {$(-1,2)$};
  \draw (0+\a+3*\a,-3*\b) node {$(0,1)$};
  \draw (0+2*\a+3*\a,\b) node {$(2,-3)$};
  \draw (0+2*\a+3*\a,0) node {$(-2,5)$};
  \draw (0+2*\a+3*\a,-\b) node {$(1,-4)$};
  \draw (0+2*\a+3*\a,2*\b) node {$(-1,4)$};
  \draw (0+2*\a+3*\a,3*\b) node {$(2,-5)$};
  \draw (0+2*\a+3*\a,4*\b) node {$(-2,3)$};
  \draw (0+2*\a+3*\a,-2*\b) node {$(-1,2)$};
  \draw (0+2*\a+3*\a,-3*\b) node {$(1,-2)$};
  \draw (0+2*\a+3*\a,-4*\b) node {$(0,1)$};

  \draw (0+0.7*\c+3*\a,0) -- (0+\a-\d+3*\a,0);
  \draw (0+0.2*\a+3*\a,0+0.2*3*\b) -- (0+0.8*\a+3*\a,0+0.8*3*\b);
  \draw (0+0.2*\a+3*\a,0-0.2*3*\b) -- (0+0.8*\a+3*\a,0-0.8*3*\b);

  \draw (\a+0.3*\a+3*\a,3*\b+0.3*\b) -- (2*\a-0.3*\a+3*\a,4*\b-0.3*\b);
  \draw (\a+0.6\a+3*\a,3*\b) -- (2*\a-0.3*\a+3*\a,3*\b);
  \draw (\a+0.3*\a+3*\a,3*\b-0.3*\b) -- (2*\a-0.3*\a+3*\a,2*\b+0.3*\b);

  \draw (\a+0.3*\a+3*\a,-3*\b-0.3*\b) -- (2*\a-0.3*\a+3*\a,-4*\b+0.3*\b);
  \draw (\a+0.6\a+3*\a,-3*\b) -- (2*\a-0.3*\a+3*\a,-3*\b);
  \draw (\a+0.3*\a+3*\a,-3*\b+0.3*\b) -- (2*\a-0.3*\a+3*\a,-2*\b-0.3*\b);

  \draw (\a+0.3*\a+3*\a,0+0.3*\b) -- (2*\a-0.3*\a+3*\a,\b-0.3*\b);
  \draw (\a+0.3*\a+3*\a,0) -- (2*\a-0.3*\a+3*\a,0);
  \draw (\a+0.3*\a+3*\a,0-0.3*\b) -- (2*\a-0.3*\a+3*\a,-\b+0.3*\b);

\end{tikzpicture}

\noindent Comparing the above two trees shows that the second tree yields the B\'{e}zout coefficients for entries in the first tree. This is not completely true for the B\'{e}zout tree of $(2,1)$ generated by $(0,1)$ (which is the same as the second tree in the above example). To fix the situation, simply change the top entry in the second level (i.e. at depth $1$) from $(-1,2)$ into $(1,-1)$, then propagate accordingly using the functions defined by $g_i(u,v),i=1,2,3$. We make this precise by introducing the following trees up to a depth of $2$ (cf. Figure 2.1 of [2], where there are typos regarding two entries in the upper right subtree).\\

\noindent {\bf Example 1.3.}\\

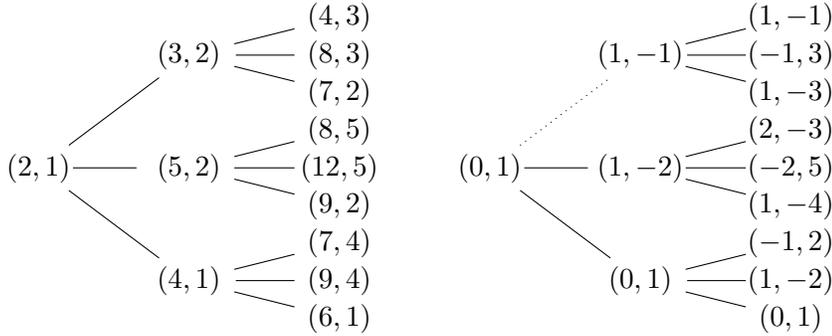
\noindent \begin{tikzpicture}

  \draw (0,0) node {$(2,1)$};
  \draw (0+1*\a,0) node {$(5,2)$};
  \draw (0+\a,3*\b) node {$(3,2)$};
  \draw (0+\a,-3*\b) node {$(4,1)$};
  \draw (0+2*\a,\b) node {$(8,5)$};
  \draw (0+2*\a,0) node {$(12,5)$};
  \draw (0+2*\a,-\b) node {$(9,2)$};
  \draw (0+2*\a,2*\b) node {$(7,2)$};
  \draw (0+2*\a,3*\b) node {$(8,3)$};
  \draw (0+2*\a,4*\b) node {$(4,3)$};
  \draw (0+2*\a,-2*\b) node {$(7,4)$};
  \draw (0+2*\a,-3*\b) node {$(9,4)$};
  \draw (0+2*\a,-4*\b) node {$(6,1)$};

  \draw (0+0.7*\c,0) -- (0+\a-\d,0);
  \draw (0+0.2*\a,0+0.2*3*\b) -- (0+0.8*\a,0+0.8*3*\b);
  \draw (0+0.2*\a,0-0.2*3*\b) -- (0+0.8*\a,0-0.8*3*\b);

  \draw (\a+0.3*\a,3*\b+0.3*\b) -- (2*\a-0.3*\a,4*\b-0.3*\b);
  \draw (\a+0.6\a,3*\b) -- (2*\a-0.3*\a,3*\b);
  \draw (\a+0.3*\a,3*\b-0.3*\b) -- (2*\a-0.3*\a,2*\b+0.3*\b);

  \draw (\a+0.3*\a,-3*\b-0.3*\b) -- (2*\a-0.3*\a,-4*\b+0.3*\b);
  \draw (\a+0.6\a,-3*\b) -- (2*\a-0.3*\a,-3*\b);
  \draw (\a+0.3*\a,-3*\b+0.3*\b) -- (2*\a-0.3*\a,-2*\b-0.3*\b);

  \draw (\a+0.3*\a,0+0.3*\b) -- (2*\a-0.3*\a,\b-0.3*\b);
  \draw (\a+0.3*\a,0) -- (2*\a-0.3*\a,0);
  \draw (\a+0.3*\a,0-0.3*\b) -- (2*\a-0.3*\a,-\b+0.3*\b);

  \draw (0+3*\a,0) node {$(0,1)$};
  \draw (0+1*\a+3*\a,0) node {$(1,-2)$};
  \draw (0+\a+3*\a,3*\b) node {$(1,-1)$};
  \draw (0+\a+3*\a,-3*\b) node {$(0,1)$};
  \draw (0+2*\a+3*\a,\b) node {$(2,-3)$};
  \draw (0+2*\a+3*\a,0) node {$(-2,5)$};
  \draw (0+2*\a+3*\a,-\b) node {$(1,-4)$};
  \draw (0+2*\a+3*\a,2*\b) node {$(1,-3)$};
  \draw (0+2*\a+3*\a,3*\b) node {$(-1,3)$};
  \draw (0+2*\a+3*\a,4*\b) node {$(1,-1)$};
  \draw (0+2*\a+3*\a,-2*\b) node {$(-1,2)$};
  \draw (0+2*\a+3*\a,-3*\b) node {$(1,-2)$};
  \draw (0+2*\a+3*\a,-4*\b) node {$(0,1)$};

  \draw (0+0.7*\c+3*\a,0) -- (0+\a-\d+3*\a,0);
  \draw[dotted](0+0.2*\a+3*\a,0+0.2*3*\b) -- (0+0.8*\a+3*\a,0+0.8*3*\b);
  \draw (0+0.2*\a+3*\a,0-0.2*3*\b) -- (0+0.8*\a+3*\a,0-0.8*3*\b);

  \draw (\a+0.3*\a+3*\a,3*\b+0.3*\b) -- (2*\a-0.3*\a+3*\a,4*\b-0.3*\b);
  \draw (\a+0.6\a+3*\a,3*\b) -- (2*\a-0.3*\a+3*\a,3*\b);
  \draw (\a+0.3*\a+3*\a,3*\b-0.3*\b) -- (2*\a-0.3*\a+3*\a,2*\b+0.3*\b);

  \draw (\a+0.3*\a+3*\a,-3*\b-0.3*\b) -- (2*\a-0.3*\a+3*\a,-4*\b+0.3*\b);
  \draw (\a+0.6\a+3*\a,-3*\b) -- (2*\a-0.3*\a+3*\a,-3*\b);
  \draw (\a+0.3*\a+3*\a,-3*\b+0.3*\b) -- (2*\a-0.3*\a+3*\a,-2*\b-0.3*\b);

  \draw (\a+0.3*\a+3*\a,0+0.3*\b) -- (2*\a-0.3*\a+3*\a,\b-0.3*\b);
  \draw (\a+0.3*\a+3*\a,0) -- (2*\a-0.3*\a+3*\a,0);
  \draw (\a+0.3*\a+3*\a,0-0.3*\b) -- (2*\a-0.3*\a+3*\a,-\b+0.3*\b);

\end{tikzpicture}

\noindent The dotted line in the second tree above means that the entry $(1,-1)$ does not come from $(0,1)$ by applying the function $g_1(u,v)$, instead one defines the entry $(1,-1)$ using the B\'{e}zout coefficients of the corresponding Pythagorean pair $(3,2)$. After this modification, it appears that the new B\'{e}zout tree yields all the B\'{e}zout coefficients of the first tree. Hence the merit of this construction is that (if it is proven to be true) it gives an efficient way to construct the B\'{e}zout coefficients for all Pythagorean pairs. \\

\noindent We can state Conjecture 2.1 of [2] (and now a theorem) as follows.\\

\noindent {\bf Theorem 1.4.} Consider the trinary trees generated by $(2,1)$ and $(3,1)$. Let $(u,v)$ be the pair in the B\'{e}zout tree corresponding to the relatively prime pair $(m,n)$ and $(U,V)$ be the pair given by the \verb"gcd" function for the same pair $(m,n)$. Then the following hold:\\
\indent (1) For all $(u,v)$ in the B\'{e}zout tree of $(3,1)$ generated by $(0,1)$, $(u,v)=(U,V).$\\
\indent (2) One third of the $(u,v)$ in the B\'{e}zout tree of $(2,1)$ generated by $(0,1)$ are not equal to $(U,V).$ Changing the value of $g(0,-1)$ in the second level of the B\'{e}zout tree from $(-1,2)$ to $(1,-1)$ results in a tree in which $(u,v)=(U,V)$ for all $(u,v).$\\

\noindent The above theorem is clearly implied by the following theorem.\\

\noindent {\bf Theorem 1.5.} For a relatively prime Pythagorean pair $(m,n)$ with $m>n>0$, except for $(m,n)=(2,1)$ and for $f_1(m,n):=f(m,-n)$, the following diagram is commutative, i.e. $\beta(f_i(m,n))=g_i(\beta(m,n)),i=1,2,3$:

\begin{center}
\begin{tikzpicture}{c}[scale=1.5]
\node (A) at (0,1) {$(m,n)$};
\node (B) at (2,1) {$(m',n')$};
\node (C) at (0,0) {$(r,s)$};
\node (D) at (2,0) {$(r',s')$};
\path[->,font=\scriptsize]
(A) edge node[above]{$f_i$} (B)
(A) edge node[right]{$\beta$} (C)
(B) edge node[right]{$\beta$} (D)
(C) edge node[above]{$g_i$} (D);
\end{tikzpicture}
\end{center}

\noindent where $\beta(m,n)$ gives the B\'{e}zout coefficients $(r,s)$ for $(m,n)$.\\

\noindent We will prove Theorem 1.5 in Section 2 using the standard Euclidean algorithm. In Section 3, we give a generalization (see Theorem 3.2).\\

\noindent {\bf 2. Euclidean Algorithm and the Proof of Theorem 1.5} \\

\noindent For simplicity we will assume that all ordered pairs $(m,n)$ consist of relatively prime integers, even though the result can be generalized to the case when $\gcd(m,n)=d>1$.\\

\noindent {\bf 2.1 Euclidean Algorithm}\\

\noindent We recall that for relatively prime integers $m>n>0$, the Division Algorithm is given by $$m=q_1n+r_1$$ $$n=q_2r_1+r_2$$ $$\cdots$$ $$r_{k-2}=q_k\cdot r_{k-1}+r_k,$$ where $r_{k-1}=\gcd(m,n)=1$ and $r_k=0$. We can record the process using matrices as follows: $$\left[\begin{array}{c}m\\ n\end{array}\right]=\left[\begin{array}{cc}q_1&1\\ 1&0\end{array}\right]\left[\begin{array}{c}n\\ r_1\end{array}\right]=\left[\begin{array}{cc}q_1&1\\ 1&0\end{array}\right]\left[\begin{array}{cc}q_2&1\\ 1&0\end{array}\right]\left[\begin{array}{c}r_1\\ r_2\end{array}\right]$$
$$=\left[\begin{array}{cc}q_1&1\\ 1&0\end{array}\right]\cdots\left[\begin{array}{cc}q_k&1\\ 1&0\end{array}\right]\left[\begin{array}{c}1\\ 0\end{array}\right].$$ Similarly for $n>m>0$ and relatively prime, we have $$\left[\begin{array}{c}m\\ n\end{array}\right]=\left[\begin{array}{cc}0&1\\ 1&0\end{array}\right]\left[\begin{array}{cc}q_1&1\\ 1&0\end{array}\right]\cdots\left[\begin{array}{cc}q_k&1\\ 1&0\end{array}\right]\left[\begin{array}{c}1\\ 0\end{array}\right].$$ Note that these intermediate matrices with left upper corner entry $q_i$ or $0$ are uniquely determined.\\

\noindent To prove Theorem 1.5, we need a few lemmas.\\

\noindent {\bf Lemma 2.2} Assume that $m>n>0$ and $\gcd(m,n)=1$. Then $(r,s)$ gives the B\'{e}zout coefficients for $(m,n)$, i.e. $(r,s)=\beta(m,n)$ if and only if $(s,r-2s)=\beta(2m+n,m).$\\

\noindent {\it Proof.} By division algorithm, we can write $$\left[\begin{array}{c}m\\ n\end{array}\right]=A\left[\begin{array}{c}1\\ 0\end{array}\right],$$ where $A$ is an invertible integer matrix of the form $$A=\left[\begin{array}{cc}q_1&1\\ 1&0\end{array}\right]\cdots \left[\begin{array}{cc}q_k&1\\ 1&0\end{array}\right],~\det(A)=\pm 1.$$ By the given assumption, $A^{-1}$ is of the form $$A^{-1}=\left[\begin{array}{cc}r&s\\ *&*\end{array}\right].$$ Now we perform the division algorithm for $2m+n$ and $m$, where the first step is the following: $$2m+n=2\cdot m+n,$$ which is followed by the division of $m$ by $n$. Hence we can write $$\left[\begin{array}{c}2m+n\\ m\end{array}\right]=\left[\begin{array}{cc}2&1\\ 1&0\end{array}\right]A\left[\begin{array}{c}1\\ 0\end{array}\right],$$ where the first row of the matrix $\left(\left[\begin{array}{cc}2&1\\ 1&0\end{array}\right]A\right)^{-1}$ gives the B\'{e}zout coefficient $(r',s')$ of the division of $2m+n$ by $m$. But $$\left(\left[\begin{array}{cc}2&1\\ 1&0\end{array}\right]A\right)^{-1}=A^{-1}\left[\begin{array}{cc}2&1\\ 1&0\end{array}\right]^{-1}$$ $$=\left[\begin{array}{cc}r&s\\ *&*\end{array}\right] \left[\begin{array}{cc}0&1\\ 1&-2\end{array}\right]=\left[\begin{array}{cc}s&r-2s\\ *&*\end{array}\right],$$ from which we have $(r',s')=(s,r-2s)$ as required. $\Box$\\

\noindent {\bf Lemma 2.3} Assume that $m>n>0$ and $\gcd(m,n)=1$. Then $(r,s)=\beta(m,m-n)$ if and only if $(s,r-s)=\beta(2m-n,m).$\\

\noindent {\it Proof.} Clearly we have $\gcd(m,m-n)=\gcd(2m-n,m)=1$. Similar to the proof of Lemma 2.2, we may write $$\left[\begin{array}{c}m\\ m-n\end{array}\right]=A\left[\begin{array}{c}1\\ 0\end{array}\right],$$ where the first row of $A^{-1}$ gives the B\'{e}zout coefficients for the division of $m$ by $m-n$, i.e. $$A^{-1}=\left[\begin{array}{cc}r&s\\ *&*\end{array}\right].$$ Performing the first step of the division of $2m-n$ by $m$, we have $$2m-n=1\cdot m+(m-n),$$ i.e. $$\left[\begin{array}{c}2m-n\\ m\end{array}\right]=\left[\begin{array}{cc}1&1\\ 1&0\end{array}\right]\left[\begin{array}{c}m\\ m-n\end{array}\right]=\left[\begin{array}{cc}1&1\\ 1&0\end{array}\right]A\left[\begin{array}{c}1\\ 0\end{array}\right],$$ whence the B\'{e}zout coefficients for the division of $2m-n$ by $m$ is given by the first row of the matrix $$\left(\left[\begin{array}{cc}1&1\\ 1&0\end{array}\right]A\right)^{-1}=A^{-1}\left[\begin{array}{cc}1&1\\ 1&0\end{array}\right]^{-1}$$
$$=\left[\begin{array}{cc}r&s\\ *&*\end{array}\right]\left[\begin{array}{cc}0&1\\ 1&-1\end{array}\right]=\left[\begin{array}{cc}s&r-s\\ *&*\end{array}\right],$$ as required. $\Box$\\

\noindent {\bf Lemma 2.4.} Assume that $m>n>0$, $n<\frac m 2$ and $\gcd(m,n)=1$. Then $$\left[\begin{array}{c}m\\ n\end{array}\right]=\left[\begin{array}{cc}q_1&1\\ 1&0\end{array}\right]A\left[\begin{array}{c}1\\ 0\end{array}\right]$$ and
$$\left[\begin{array}{c}m\\ m-n\end{array}\right]=\left[\begin{array}{cc}1&1\\ 1&0\end{array}\right]\left[\begin{array}{cc}q_2'&1\\ 1&0\end{array}\right]A\left[\begin{array}{c}1\\ 0\end{array}\right],$$ where $q_1>1$ and $1+q_2'=q_1$.\\

\noindent {\it Proof.} By assumption, we can write $$m=q_1n+r$$ with $q_1=\lfloor\frac m n\rfloor\geq 2$ and $r<n$.\\

\noindent Letting $q_2'=q_1-1$, we have $$m=1\cdot (m-n)+n,$$ where $n<m-n$ by assumption, and $$m-n=(q_1-1)n+r=q_2'n+r.$$ Writing $$\left[\begin{array}{c}n\\ r\end{array}\right]=A\left[\begin{array}{c}1\\ 0\end{array}\right]$$ and expressing the above divisions in terms of matrices, the result is clear. $\Box$\\

\noindent {\bf Lemma 2.5} Assume that $m>n>0$, $\gcd(m,n)=1$ and $(m,n)\neq (2,1)$. Then $(r,s)=\beta(m,n)$ if and only if $(r+s,-s)=\beta(m,m-n)$.\\

\noindent {\it Proof.} Since $(m,n)\neq (2,1)$, there are only the following two cases to consider.\\
\noindent {\it Case 1:} $n<\frac m 2$. Using Lemma 2.4, we see that the division of $m$ by $n$ is described by the procedure $$\left[\begin{array}{cc}q_1&1\\ 1&0\end{array}\right]A,\eqno (1)$$ if and only if the division of $m$ by $m-n$ is described by the following procedure: $$\left[\begin{array}{cc}1&1\\ 1&0\end{array}\right]\left[\begin{array}{cc}q_2'&1\\ 1&0\end{array}\right]A,\eqno (2)$$ where $$q_1=q_2'+1.$$ The fact that $(r,s)$ is the B\'{e}zout coefficient for division of $m$ by $n$ means precisely that $$A^{-1}\left[\begin{array}{cc}q_1&1\\ 1&0\end{array}\right]^{-1}=\left[\begin{array}{cc}r&s\\ *&*\end{array}\right],$$ while the B\'{e}zout coefficient for the division of $m$ by $m-n$ is given by the first row vector of the matrix $$A^{-1}\left[\begin{array}{cc}q_2'&1\\ 1&0\end{array}\right]^{-1}\left[\begin{array}{cc}1&1\\ 1&0\end{array}\right]^{-1}.$$ But
$$A^{-1}\left[\begin{array}{cc}q_2'&1\\ 1&0\end{array}\right]^{-1}\left[\begin{array}{cc}1&1\\ 1&0\end{array}\right]^{-1}$$
$$=A^{-1}\left[\begin{array}{cc}q_1-1&1\\ 1&0\end{array}\right]^{-1}\left[\begin{array}{cc}1&1\\ 1&0\end{array}\right]^{-1}$$
$$=A^{-1}\left(\left[\begin{array}{cc}1&-1\\ 0&1\end{array}\right]\left[\begin{array}{cc}q_1&1\\ 1&0\end{array}\right]\right)^{-1}\left[\begin{array}{cc}0&1\\ 1&-1\end{array}\right]$$
$$=A^{-1}\left[\begin{array}{cc}q_1&1\\ 1&0\end{array}\right]^{-1}\left[\begin{array}{cc}1&1\\ 0&1\end{array}\right]\left[\begin{array}{cc}0&1\\ 1&-1\end{array}\right]$$
$$=\left[\begin{array}{cc}r&s\\ *&*\end{array}\right]\left[\begin{array}{cc}1&0\\ 1&-1\end{array}\right]=\left[\begin{array}{cc}r+s&-s\\ *&*\end{array}\right],$$ which shows the result. \\

\noindent \emph{Case 2:} $n>\frac m 2$. Assume that $(r,s)=\beta(m,n)$ and $(r',s')=\beta(m,m-n).$ Since $m-n<\frac m 2,$ by the result of Case 1, we have $$r=r'+s'~ {\rm and~}s=-s',$$ which shows that $r'=r+s$ and $s'=-s$, as required. $\Box$\\

\noindent {\bf Remark 2.6} The above lemma does not hold when $(m,n)=(2,1)$, since $(0,1)=\beta(m,n)=\beta(m,m-n)$ here, but $(r,s)=(0,1)\neq (r+s,-s).$\\

\noindent {\bf Lemma 2.7} Assume that $m>n>0, \gcd(m,n)=1$ and $(m,n)\neq (2,1)$. Then $(r,s)=\beta(m,n)$ if and only if $(-s,r+2s)=\beta(2m-n,m).$\\

\noindent {\it Proof.} The mapping $(m,n)\mapsto (2m-n,m)$ can be factored as $(m,n)\mapsto (m,m-n)\mapsto (2m-n,m)$, so by Lemma 2.5 and Lemma 2.3, the corresponding B\'{e}zout coefficients are given by $(r,s)\mapsto (r+s,-s)\mapsto (-s,(r+s)-(-s))=(-s,r+2s)$, as required. $\Box$\\

\noindent {\bf Proof of Theorem 1.5}\\

\noindent The proof for the pair $f_1$ and $g_1$ follows from Lemma 2.7. Note that the condition $(m,n)\neq (2,1)$ is precisely used here. The proof for the pair $f_2$ and $g_2$ follows from Lemma 2.2. The proof for the pair $f_3$ and $g_3$ is essentially the same as that of the previous case. Here are the details. Let $$\left[\begin{array}{c}m\\ n\end{array}\right]=\left[\begin{array}{cc}q_1&1\\ 1&0\end{array}\right]A\left[\begin{array}{c}1\\ 0\end{array}\right]$$ where the first step of the division process is written out. Now for the division of $2n+m$ by $n$, one has $$\left[\begin{array}{c}2n+m\\ n\end{array}\right]=\left[\begin{array}{cc}q_1+2&1\\ 1&0\end{array}\right]A\left[\begin{array}{c}1\\ 0\end{array}\right].$$ Let $(r,s)$ be the B\'{e}zout coefficients for the division of $m$ by $n$, i.e. $$\left[\begin{array}{cc}r&s\\ *&*\end{array}\right]=A^{-1}\left[\begin{array}{cc}q_1&1\\ 1&0\end{array}\right]^{-1}.$$ If $(r',s')$ is the B\'{e}zout coefficients for the division of $2n+m$ by $n$, then $$\left[\begin{array}{cc}r'&s'\\ *&*\end{array}\right]=A^{-1}\left[\begin{array}{cc}q_1+2&1\\ 1&0\end{array}\right]^{-1}$$
$$=A^{-1}\left(\left[\begin{array}{cc}1&2\\ 0&1\end{array}\right]\left[\begin{array}{cc}q_1&1\\ 1&0\end{array}\right]\right)^{-1}$$
$$=A^{-1}\left[\begin{array}{cc}q_1&1\\ 1&0\end{array}\right]^{-1}\left[\begin{array}{cc}1&2\\ 0&1\end{array}\right]^{-1}$$
$$=\left[\begin{array}{cc}r&s\\ *&*\end{array}\right]\left[\begin{array}{cc}1&-2\\ 0&1\end{array}\right]=\left[\begin{array}{cc}r&s-2r\\ *&*\end{array}\right],$$ so $$r'=r,s'=s-2r,$$ as required. $\Box$\\

\noindent {\bf 3. A Generalization}\\

\noindent We first extend the definition of B\'{e}zout coefficients to general ordered pairs $(m,n)\in {\mathbb Z}\times {\mathbb Z}$. The following definition seems to yield the same output as the \textsc{Matlab}'s function $[G,U,V]={\verb"gcd"}(m,n)$ [3] or \textsc{Sage}'s ${\verb"xgcd"}$ function [6]. We have tested this by writing a Sage script using the following definitions for relatively prime $(m,n)$ up to a reasonable size. In any case, our proof will be based on the following definitions.\\

\noindent {\bf Definition 3.1.} The B\'{e}zout coefficients $\beta(a,b)$ for an ordered pair $(a,b)\in {\mathbb Z}\times{\mathbb Z}$ are defined as follows:\\

\noindent {\bf 3.1.1} For $a>b>0$, $\beta(a,b)=(r,s)$ with $ra+sb=\gcd(a,b)$ is given by the Euclidean algorithm which is uniquely determined. One writes this as $(a,b)\mapsto (\gcd(a,b),r,s)$. \\

This is extended to all ordered pairs by the following rules:\\

\noindent {\bf 3.1.2} $(0,a)\mapsto (|a|,0,{\rm sign}(a))$, where ${\rm sign}(a)$ is the sign of $a$, and by convention ${\rm sign}(0)=0$\\

\noindent {\bf 3.1.3} $(\pm a,a)\mapsto (|a|,0,{\rm sign}(a))$\\

\noindent {\bf 3.1.4} If $|a|\neq |b|$ and $\beta(|a|,|b|)=(r,s),$ then $(a,b)\mapsto (\gcd(a,b),{\rm sign}(a)r,{\rm sign}(b)s).$\\

\noindent {\bf 3.1.5} If $|a|\neq |b|$ and $(a,b)\mapsto (\gcd(a,b),r,s)$, then $(b,a)\mapsto (\gcd(a,b),s,r)$\\

\noindent We leave the readers to check that these formulas are consistent.\\

\noindent {\bf Theorem 3.2} Let $A$ be a unimodular $2\times 2$ matrix, i.e. an integer matrix such that $\det(A)=\pm 1$. Consider the mappings $$f\left(\left[\begin{array}{c}m\\ n\end{array}\right]\right):=A\left[\begin{array}{c}m\\ n\end{array}\right]$$ and $$g\left(\left[\begin{array}{c}m\\ n\end{array}\right]\right):=(A^{-1})^T\left[\begin{array}{c}m\\ n\end{array}\right]$$ for $(m,n)\in {\mathbb Z}\times {\mathbb Z}$ such that $\gcd(m,n)=1$. Then with only finitely many exceptions of relatively prime ordered pairs $(m,n)$, one has $g(\beta(m,n))=\beta(f(m,n))$, where $\beta$ maps an ordered pair $(m,n)\in {\mathbb Z}\times {\mathbb Z}$ to its B\'{e}zout coefficients.\\

\noindent {\it Proof.} The idea of the proof is that if the group of unimodular $2\times 2$ matrices, denoted ${\rm GL}_2({\mathbb Z})$, is finitely generated, then we can decompose each element $A$ in the group as a product of its generators, say $$A=V_kV_{k-1}\cdot V_1,$$ where each $V_i$ lies in a finite set of generators. Since $$(A^{-1})^T=(V_k^{-1})^T(V_{k-1}^{-1})^T\cdots (V_1^{-1})^T,$$ the proof of compatibility of B\'{e}zout coefficients of relatively prime ordered pairs under the transformation $A$ is reduced to the simple case when $A=V_i$, where $V_i$ is in the set of generators, and we check the relation $\beta(f(m,n))=g(\beta(m,n))$ for $f$ defined by $V_i$ and for $g$ defined by $(V_i^{-1})^T$. This is because for a factorization $A$ into the generators (so $A$ is a series of compositions of the generator functions), if compatibility holds at each step of the successive composition with a finite number of exceptions, then it is easy to see that there will be only finitely many exceptions for the final composite function, which is $A$ (we illustrate this in Example 3.3). Now we start to prove the result for its generators.\\

\noindent It is well known [1] that ${\rm GL}_2({\mathbb Z})$ is generated by $${\rm U}=\left[\begin{array}{cc}0&1\\ 1&0\end{array}\right],{\rm S}=\left[\begin{array}{cc}0&-1\\ 1&0\end{array}\right]~{\rm and~} {\rm T}=\left[\begin{array}{cc}1&1\\ 0&1\end{array}\right].$$ It is easy to check that the sets of exceptional pairs for ${\rm U},{\rm S}~{\rm and~}{\rm S}^{-1}$ are all given by $$\{(-1,-1),(-1,1),(1,-1),(1,1)\}.$$

\noindent To determine the exceptional set $E_{\rm T}$ for ${\rm T}$, we check first the following special cases $(m,n)$ such that
$$(m,n)\in \{(0,\pm 1),(\pm 1,0),(\pm 1,\pm 1),(\pm 1,\mp 1)\}$$ or $$(m+n,n)\in \{(0,\pm 1),(\pm 1,0),(\pm 1,\pm 1),(\pm 1,\mp 1)\}$$ This gives exceptional ordered pairs $(\pm 1,0)$, for which $g(\beta(m,n))\neq \beta(f(m,n))$.\\

\noindent For the remaining cases, we may assume that $$(|m|,|n|),(|m+n|,|n|)\notin\{(1,0),(0,1),(1,1)\}.$$ After excluding the special cases above, we use 3.14 and 3.15 to reduce the checking to the following cases, noting that $\beta(-m,-n)=-\beta(m,n)$ and $\beta(-m-n,-n)=-\beta(m+n,n)$. \\

\noindent {\it Case:} $n>m>0$. Let $$\left[\begin{array}{c}m\\ n\end{array}\right]=\left[\begin{array}{cc}0&1\\ 1&0\end{array}\right]A\left[\begin{array}{c}1\\ 0\end{array}\right],\left[\begin{array}{c}m+n\\ n\end{array}\right]=\left[\begin{array}{cc}1&1\\ 1&0\end{array}\right]A\left[\begin{array}{c}1\\ 0\end{array}\right].$$ Then $$\left[\begin{array}{cc}r'&s'\\ *&*\end{array}\right]=A^{-1}\left[\begin{array}{cc}1&1\\ 1&0\end{array}\right]^{-1}=A^{-1}\left[\begin{array}{cc}0&1\\ 1&0\end{array}\right]\left[\begin{array}{cc}1&-1\\0&1\end{array}\right]$$
$$=\left[\begin{array}{cc}r&s\\ *&*\end{array}\right]\left[\begin{array}{cc}1&-1\\0&1\end{array}\right]=\left[\begin{array}{cc}r&s-r\\ *&*\end{array}\right],$$ where $(r,s)=\beta(m,n)$ and $(r',s')=\beta(m+n,n)$.\\

\noindent {\it Case:} $m>n>0$. \\

\noindent Let $$\left[\begin{array}{c}m\\ n\end{array}\right]=\left[\begin{array}{cc}q_1&1\\ 1&0\end{array}\right]A\left[\begin{array}{c}1\\ 0\end{array}\right].$$ Then $$\left[\begin{array}{cc}r&s\\ *&*\end{array}\right]=A^{-1}\left[\begin{array}{cc}q_1&1\\ 1&0\end{array}\right]^{-1},$$ $$\left[\begin{array}{c}m+n\\ n\end{array}\right]=\left[\begin{array}{cc}q_1+1&1\\ 1&0\end{array}\right]A\left[\begin{array}{c}1\\ 0\end{array}\right],$$ hence $$\left[\begin{array}{cc}r'&s'\\ *&*\end{array}\right]=A^{-1}\left[\begin{array}{cc}q_1+1&1\\ 1&0\end{array}\right]^{-1}$$
$$=A^{-1}\left[\begin{array}{cc}q_1&1\\ 1&0\end{array}\right]^{-1}\left[\begin{array}{cc}1&1\\ 0&1\end{array}\right]^{-1}$$
$$=\left[\begin{array}{cc}r&s\\ *&*\end{array}\right]\left[\begin{array}{cc}1&-1\\ 0&1\end{array}\right]=\left[\begin{array}{cc}r&s-r\\ *&*\end{array}\right].$$ It follows that for both of the above cases, $$\left[\begin{array}{c}r'\\ s'\end{array}\right]=\left(\left[\begin{array}{cc}1&1\\ 0&1\end{array}\right]^{-1}\right)^T\left[\begin{array}{c}r\\ s\end{array}\right],$$ as required.\\

\noindent {\it Case:} $m>0, n<0.$ Let $n=-n'$. \\

\noindent {\it Subcase 1:} $m<n'$.\\

\noindent Tracing the relations $$\left[\begin{array}{c}m\\ n\end{array}\right]\leftrightarrow \left[\begin{array}{c}m\\ n'\end{array}\right]\leftrightarrow \left[\begin{array}{c}n'\\ m\end{array}\right]$$ and $$\left[\begin{array}{c}m+n\\ n\end{array}\right]\leftrightarrow \left[\begin{array}{c}n'-m\\ n'\end{array}\right]\leftrightarrow \left[\begin{array}{c}n'\\ n'-m\end{array}\right],$$ it suffices to find the relation $$\left[\begin{array}{c}n'\\ m\end{array}\right]\leftrightarrow \left[\begin{array}{c}n'\\ n'-m\end{array}\right].$$ When $(n',m)\neq (2,1)$ (i.e. when $(m,n)\neq (1,-2)$), this can be determined by Lemma 2.5. Using 3.1.4, 3.1.5 and Lemma 2.5, we find the same relation between $(r',s')$ and $(r,s)$ as above. The special cases $(m,n)=(\pm 1,\mp 2)$ are checked directly to be exceptional.\\

\noindent {\it Subcase 2:} $m>n'$.\\

\noindent Tracing the relations $$\left[\begin{array}{c}m\\ n\end{array}\right]\leftrightarrow \left[\begin{array}{c}m\\ n'\end{array}\right]$$ and $$\left[\begin{array}{c}m+n\\ n\end{array}\right]\leftrightarrow \left[\begin{array}{c}m-n'\\ n'\end{array}\right],$$ it suffices to find the relation $$\left[\begin{array}{c}m\\ n'\end{array}\right]\leftrightarrow \left[\begin{array}{c}m-n'\\ n'\end{array}\right].$$ We have $$\left[\begin{array}{c}m\\ n'\end{array}\right]=\left[\begin{array}{cc}q_1'+1&1\\ 1&0\end{array}\right]A\left[\begin{array}{c}1\\ 0\end{array}\right]$$ and $$\left[\begin{array}{c}m-n'\\ n'\end{array}\right]=\left[\begin{array}{cc}q_1'&1\\ 1&0\end{array}\right]A\left[\begin{array}{c}1\\ 0\end{array}\right],$$ where $q_1'=0$ if $m<2n'$ and $q_1'\geq 1$ if $m\geq 2n'$. As a result, we find the same relation as above. In summary, for the transformation ${\rm T}$, the set $E_{\rm T}$ of exceptional cases is given by $$E_{\rm T}=\{(-1,0),(1,0),(1,-2),(-1,2)\}.$$ Similarly the exceptional set $E_{{\rm T}^{-1}}$ is given by $$E_{{\rm T}^{-1}}=\{(-1,-2),(-1,0),(1,0),(1,2)\}.$$ This concludes the proof. $\Box$\\

\noindent {\bf Example 3.3} The factorization of $$\left[\begin{array}{cc}1&2\\ 0&1\end{array}\right]={\rm T}^2$$ $$\left[\begin{array}{cc}2&1\\ 1&0\end{array}\right]={\rm T}^2{\rm U}$$ and $$\left[\begin{array}{cc}2&-1\\ 1&0\end{array}\right]={\rm T}^2{\rm S}$$ allows us to determine the exceptional set of these transformations. For example, using the proof of the above theorem, let's determine the exceptional set $E_{{\rm T}^2{\rm S}}$ for the transformation $\left[\begin{array}{cc}2&-1\\ 1&0\end{array}\right]$, which is described by the following process:
$$\begin{array}{ccccccc} (m,n)&\rightarrow&{\rm S}(m,n)&\rightarrow&{\rm T}({\rm S}(m,n))&\rightarrow&{\rm T}({\rm TS}(m,n))\\
                           {\rm (a)}&\rightarrow&{\rm (b)}&         \rightarrow&{\rm (c) }                 &\rightarrow&{\rm (d)}\end{array},$$ where compatibility of B\'{e}zout coefficients can fail at (a) for ordered pairs in the exceptional set $E_{{\rm S}}$ of $S$, or at (b) for ordered pairs in the exceptional set $E_{\rm T}$ of ${\rm T}$, or at (c) for ordered pairs in the exceptional set $E_{{\rm T}}$ of ${\rm T}$. Taking preimage of these exceptional sets to the beginning step (a), we see that the compatibility of B\'{e}zout coefficients for ${\rm T}^2{\rm S}$ can only possibly fail for ordered pairs in the set $$E_{\rm S}\cup {\rm S}^{-1}E_{{\rm T}}\cup ({\rm TS})^{-1}E_{\rm T}.$$ By direct checking, this turns out to be all the exceptional cases, i.e. $$E_{{\rm T}^2{\rm S}}=\{(-2,-3),(-2,-1),(-1,-1),(-1,1),(0,-1),(0,1),(1,-1),(1,1),(2,1),(2,3)\}.$$ Similarly, we get that
                           $$E_{{\rm T}^2}=\{(-3,2),(-1,0),(-1,2),(1,-2),(1,0),(3,-2)\}$$ and
                           $$E_{{\rm T}^2{\rm U}}=\{(-2,1),(-2,3),(-1,-1),(-1,1),(0,-1),(0,1),(1,-1),(1,1),(2,-3),(2,-1)\}.$$

\noindent {\bf References}\\

\noindent [1] T. Apostol, Modular functions and Dirichlet series in number theory, Springer-Verlag, 1976.\\

\noindent [2] E. Gullerud and J.S. Walker, Generalized B\'{e}zout trees for Pythagorean pairs, arXiv:1803.04875v1 [math.NT] 13 Mar 2018.\\

\noindent [3] Matlab. A language and environment for technical computing. Product of MathWorks.\\

\noindent [4] I. Niven, H.S. Zuckerman, and H.L. Montgomery, An introduction to the theory of numbers, 5th edition, John Wiley \& Sons, Inc.\\

\noindent [5] T. Randall and R. Saunders, ``The family tree of the Pythagorean triplets revisited." \emph{Math. Gaz.} ({\bf 78}) (482), pp. 190-193, 1994.\\

\noindent [6] William A. Stein et al. Sage Mathematics Software (Version 6.10), http://www.sagemath.org.\\

\indent \textsc{\small Department of Mathematics, Norfolk State University} \\
\indent {\small \it E-mail address:} {\small\verb"ctperng@nsu.edu"}\\

\indent \textsc{\small Department of Mathematics, Norfolk State University} \\
\indent {\small \it E-mail address:} {\small\verb"mcbrucal-hallare@nsu.edu"}\\

\end{document}